
\input amstex
\documentstyle{amsppt}
\magnification =\magstep 1
\refstyle{A}
\NoRunningHeads

\topmatter
\title Global properties of families of plane curves \endtitle
\author Robert  Treger \endauthor
\address AG{\&}N, Princeton, New Jersey \endaddress
\email tregrob{\@}aol.com \endemail
\abstract We describe degenerations of projective plane curves to
curves containing
a fixed line $l$ as a component, and show that $H^1({\overline
V}_{n,d,m}
,{\Cal O} (r))=0,\  r \in{\Bbb Z}$,  where  $V_{n,d,m}\subset {\Bbb
P}^N\ (N = n(n+3)/2)$
is the subscheme consisting
of irreducible plane curves having smooth contact of order at least
$m$ with $l$  at a fixed point $\bold p \in l$ and $d$ nodes and no
other singularities.
 \endabstract

\toc
\head 0. Introduction \endhead
\head 1. Preliminaries \endhead
\head 2. Splitting off a line:  a theorem \endhead
\head 3. Splitting off a line:  lemmas and a proposition \endhead
\head 4. Admissible schemes:  notation, definitions, and lemmas
\endhead
\head 5. A vanishing theorem for admissible schemes \endhead
 \endtoc
\endtopmatter

\document
\head 0. Introduction \endhead

An arbitrary smooth projective curve is birationally equivalent
to a plane nodal curve.  We denote  by  %
${\Bbb P}^N$ the projective space parametrizing all projective
plane curves of degree  $n\  (N = n(n+3)/2).$  Let  %
$l \subset {\Bbb P}^2$ be a fixed line, and  %
$ \bold p \in l$ a fixed point.  Let $ n$, $d$, $m$  be three integers
with
$0\leq m \leq n$  and  $1\leq d \leq (n-1)(n-2)/2.$
Let $V_{n,d,m}\subset {\Bbb P}^N$
be the (locally closed) subscheme consisting
of irreducible curves having smooth contact of order at least
$m$ with $l$  at $\bold p$ and $d$ nodes and no other singularities.
These schemes
are smooth and $\dim V_{n,d,m} = N-d-m$ \cite{H,~Sect.~2}.
We denote by $U_{m}(n,g)$  the {\it closure\/} in
${\Bbb P}^N$ of the locus of reduced plane curves of degree
$n$  and geometric genus  $g$, not containing $l$ and having contact
of order $m$ with
$l$ at $\bold p$.

   In this paper we study global properties of  ${\overline V}_{n,d,m}
\subset {\Bbb P}^N$ and related schemes.  In Theorem 2.2, we
describe
degenerations of a general member of a component of $U_{m}(n,g)$
 to a curve containing the line $l$. In Proposition 3.9, we restrict
ourselves to
${\overline V}_{n,d,m}$
  and obtain a more precise result.  Fortunately, even the case
when  d = 1  is not trivial (Lemma 3.8).  We prove both theorems by induction
on the number of nodes.
 Finally, in Theorem 5.2, we prove the vanishing of
$H^1(V,{\Cal O} (r))$,  $r \in
{\Bbb Z}$,  where  V  is an arbitrary  ${\overline V}_{n,d,m}$-
admissible scheme
(see Definition 4.5). In particular
$H^1({\overline V}_{n,d,m},
{\Cal O} (r))=0,\  r \in
{\Bbb Z}.$  The vanishing theorem is proved by induction on the
size of admissible schemes by taking suitable hyperplane sections.

 This paper suggests that instead of dealing with a single  ${\overline
V}_{n,d,0}$
 one should consider all  ${\overline V}_{n,d,m}$-admissible schemes
for  $0\leq m \leq n.$ We observe that ${\overline V}_{n,d,m}$-
admissible schemes
are often generically non-reduced, and this plays a key role in the
proof
of the vanishing theorem.

The author is grateful to Ching-Li Chai, David Eisenbud, Joe Harris,
Michael Larsen, Boris Moishezon, Zinovy Reichstein, Stephen Shatz
and Slava Shokurov for stimulating discussions.  Especially helpful
were remarks by Chai, Harris and Shokurov.

\head 1. Preliminaries \endhead

\subhead 1.1. Notation \endsubhead
Let $$\Sigma_{n,d,m} \subset {\Bbb P}^N \times {\text
{Sym}}^d({\Bbb P}^2)$$ be
the {\it closure\/} of the locus of pairs  $(E,\, \sum_{i=1}^dR_i),$
where  $E$  is an irreducible
nodal curve having smooth contact of order at least  $m$  with
$l$ at  $\bold p,$  and  $R_{1}$, \dots , $R_{d}$  are its nodes.  We
denote by
$\pi_{N}$  and  $\pi_{n,d}$ the projections of
the product to ${\Bbb P}^N$  and ${\text {Sym}}^d({\Bbb P}^2),$
respectively, and identify  $V_{n,d,m}$  with
$\pi_N^{-1}(V_{n,d,m}).$
For $1\leq d \leq (n-2)(n-3)/2,$ we set
$$
 V_{n,d,n+1}=\{ C+l \in {\Bbb P}^N \ \bigm | C \in V_{n-1,d,0}
\setminus (C_{\bold p}
\cup M_2)
 \subset {\Bbb P}^{N_1}\},
$$
where $N_{1}= (n-1)(n+2)/2$,  $C_{\bold p}$  is the divisor in ${\Bbb
P}^{N_1}$
of curves containing $\bold p ,$ and  $M_2$
is the divisor in
${\Bbb P}^{N_1}$
of curves having multiple points of intersection
with
$l$ ($M_2$  is described in \cite {DH1, I, Sect.~3}).  We
get a chain of closed subschemes in
${\Bbb P}^N$ (see Lemma 3.3 below):
$${\overline V}_{n,d,n+1} \subset {\overline V}_{n,d,n} \subset \dots
\subset
 {\overline V}_{n,d,0}.
$$

Let  $X$, $Y$, $Z$  be a coordinate system in
${\Bbb P}^2$ such that
$l = \{X = 0 \}$  and
${\bold p} =\mathbreak [0:1:0].$  Let
$$
f_{C}(X,Y,Z)
 = \sum a_{jk}X^{j}Y^{k}Z^{n-j-k}=  X(\cdots) + \sum a_{0k}Y^{k}Z^{n-
k}= 0
$$ be
an equation of a curve  $C.$  For  $m \geq 1,$  the condition  $a_{on}
=\dots= a_{on+1-m} = 0$
means that  $C$  has contact of
order at least  $m$  with $l$
at ${\bold p}$  (we do not exclude the case
$l \subset C$).

 The standard action of
$PGL(2)$  on
${\Bbb P}^2$ induces an action of
$PGL(2)$  on the parameter space
${\Bbb P}^{N}$,  and we can identify
$PGL(2)$  with the corresponding subgroup of
$PGL(N).$
We denote by
$
\varphi_{s_1s_2}$ ($\phi_{t_1t_2}$) the transformation of ${\Bbb
P}^N$
induced by the transformation
$$
[X:Y:Z] \rightarrow [X:Y:s_1Z+s_2X] \qquad ([X:Y:Z] \rightarrow
[X:t_1Y+t_2X:Z]).
$$

 Let %
$\bold G \subset PGL(N)$  denote the subgroup generated by all the
elements
$\varphi_{s_1s_2}$ and $\phi_{t_1t_2}$ ($s_1$, $s_2$, $t_1$, $t_2
\in {\Bbb C}^*).$

 We denote by  $g(E)$  the geometric genus of a reduced irreducible
curve  E.  The geometric genus of a reduced curve $E^{\prime} +
E^{\prime\prime}$
is defined inductively:  $g(E^{\prime} + E^{\prime\prime}) =
g(E^{\prime}) + g(E^{\prime\prime}) - 1.$  Clearly
 $g(E) \geq -\deg(E) + 1$ with equality if and only if $E$ is
a union of lines.  For a point $Q$ of an arbitrary reduced curve
 $E,$  we set $\delta_Q(E) = \dim_{\Bbb C} \widetilde{\Cal O}_Q
/{\Cal O}_Q,$
where ${\Cal O}_Q$ is the local ring of $E$ at $Q$ and
$\widetilde{\Cal O}_Q$
its normalization.

By a point of a scheme we mean a %
{\it closed\/}  point unless stated otherwise, and by a branch
(of a scheme) through a point we mean a %
{\it local\/}  branch.  By a component of a scheme we mean an
{\it isolated\/}  component.  We denote by  %
${\bold T}_{v}(W)$  the Zariski tangent space of  $W$  at a point
 $v.$  Finally, to a family  $W$  of curves of degree  $n$  in  %
${\Bbb P}^{2}$ one can canonically associate a
closed subscheme  $\overline W \subset
{\Bbb P}^{N}$;  we assume that  $\overline W$  derived from a
{\it family\/}  is always {\it reduced}.

We will frequently use the following basic facts.

\proclaim {1.2. Proposition {\rm {(\cite {H, Lemma 2.4})}}}
Let  $ W $ be a component of  $U_{m}(n,g).$
For $m \ge 1,$ we have:

{\rm {a)}} $\dim W = 3n + g - 1 - m;$ and

{\rm {b)}} if  $E$
is a general member of $ W,$ then  in a neighborhood of
$\bold p,$\ $E$
is a union of smooth arcs having contact of orders $a_1,$
\dots , $a_{r}$ with
$l,$\  $\sum a_i = m,$  and having minimal order of contact among
themselves;
the remaining singularities of  $E$  are all nodes, and  ${(E \cdot
l)}_Q\le 1$ if
$Q\ne{\bold p}$;  %
moreover  $\delta_p(E) = \sum b_{j}{(b_j -1)/2}$ and  $\sum b_j =
m,$
where  $b_1 = b_1(W)$ is the multiplicity of $ E$  at
$\bold p$
and $b_2= b_2(W),$ \dots , $ b_m = b_m(W)$ are the multiplicities of
$E$
 at  its infinitely near points lying over $\bold p$  in the direction of
$l.$
\endproclaim

\proclaim {1.3. Proposition} {\rm {i)}} $\Sigma_{n,d,0}$ is unibranch
everywhere
\cite {Tr}.

{\rm {ii)\ }} For  $d\le n,$\  $\Sigma_{n,d,m}$ is unibranch
everywhere.

{\rm {iii)}} $V_{n,d,0}$   is irreducible \cite {H}.
\endproclaim

\demo{Proof}  ii) We get $n(n+3)/2 - 3d - m \ge 0$ with equality if
and only if $n=d=m=5.$
Hence the map  $\pi_{n,d}: \Sigma_{n,d,m}\rightarrow
\text{Sym}^{d}(
{\Bbb P}^{2})$  is surjective and one can repeat an elementary
argument from \cite {Tr, Theorem (Case: $n >> d$)}.
\enddemo

\head  2. Splitting off a line:  a theorem \endhead
\subhead 2.1. Assigned and virtually non-existent nodes
\endsubhead
Let  F  be
a curve of degree  $n$  with  $d$  nodes.  If  $F$  is regarded as
the limit of a variable curve of degree  $n$  with  $d^{\prime} <d$
nodes,
then it is said that the  $d^{\prime}$  nodes of  $F$  that are very
near
the  $d^{\prime}$  nodes of the variable curve are
{\it assigned\/}  nodes of  $F,$  while the remaining  $d -
d^{\prime}$  nodes
of  $F$  are considered as
{\it virtually non-existent}.

 One may also consider assigned cycles on arbitrary curves.
An assigned  $d$-cycle  $b$  on a curve  $E$  is said to be
{\it connected\/} (in the sense of Severi) if there is an
{\it irreducible\/} nodal curve  $E^{\prime}$  whose  $d$  nodes
approach
 $b$  as  $E^{\prime}$  degenerates to  $E.$   Further, a branch of
${\overline V}_{n,d,m}$
 through an arbitrary reduced  $E \in {\overline V}_{n,d,m}$
determines
a  $d$-cycle of assigned singularities on  E.

 Let  $W$  be a component of $U_{m}(n,g)$  with  $m \ge 1.$
One can construct a maximal irreducible subfamily $W_{1} \subset
U_{m}(n,g+1)$
such that  $W \subset W_{1}.$  Take a general member  $F$  of
$U_{0}(n,g)$
 that is very near a general member of  $W.$  To obtain a family
of curves of genus  $g + 1,$  we regard one node of  $F$  as virtually
non-existent.  We then cut that family by the  $m$  hyperplanes
in $
{\Bbb P}^{N}$,   $a_{0n} = \cdots = a_{0n+1-m}  = 0,$  and take a
suitable irreducible
subfamily.

 \proclaim {2.2. Theorem} Let  $W$
be a component of  $U_{m}(n,g).$  Let  $W^{\prime}$
be a component of $W \cap \{a_{0n-m} = 0  \}$ whose general
member
is of the form $C + l.$ Then $C$
 is a general member of a component of  $U_{m^{\prime}}(n-
1,g^{\prime}),$
 where
$$
m^{\prime}= m - (g - g^{\prime}) - 2.
$$
Further, let ${\Cal F}(E,C{+}l) \subset W$ be an arbitrary local
irreducible
$1$-dimensional family through $C + l$  such that all members of
${\Cal F}(E,C{+}l)
\backslash C{+}l$, denoted by $E$,  are general points of $W.$ Then
$$
\sum_{Q \rightarrow {\bold p}} \delta_{Q}(E) =\delta_{\bold
p}(C{+}l), \quad
\sum_{Q \not \to l} \delta_{Q}(E) +
\sum_{Q \rightarrow l \backslash {\bold p}} \delta_{Q}(E) =
\sum_{R \notin l} \delta_{R}(C{+}l) + n - m.
$$
\endproclaim

         \demo{Proof}
We will prove the theorem by induction on  $n$
and $\Delta_W =\sum_{Q} \delta_{Q}(E),$
where  $E$  is a general member of  $W.$
The theorem is trivial if  $n\le 2$  or  $\Delta_W =0.$
We assume that $\Delta_{W} > 0.$  The induction is based on the
construction described in (2.1).  Let  $W_{1}\subset U_{m}(n,g+1)$
be any maximal irreducible subfamily such that  $W\subset W_{1}.$
Let
$E_{1} \in W_{1}$ be a general member
that is very near  $E.$

     A subvariety of  $W$  consisting of nonreduced curves has
codimension
strictly greater than  $1$  in  $W.$  In case general members of
the subvariety do not contain $l,$  this follows from the semi-stable
reduction theorem
for families of curves and standard dimension counts (\cite {DH1,
Sect.~1(a)}, \cite{H}, \cite{N}).  Otherwise, we take  $W_{1}\supset
W,$  as above,
 and use the induction on $n$  and  $\Delta_W.$
In particular  $C + l$ is a reduced curve.

 Since  $g(C{+}l) = g(C) - 1$  and  $\sum_{Q} \delta_{Q}(E)  \le
\sum_{Q} \delta_{Q}(C{+}l),$  we have  $g(C) = g(E) - k$  with  $k\ge -
1$ \cite {Hi}.  %

{\it Claim\/} 1:\ \ $m(C):=(C\cdot l)_{\bold p } \le m - k - 2.$
  Indeed, since  $W \nsubseteq \{a_{0n-m}
= 0 \}$,   $C$  is moving in a family of dimension  $3n + g(E) -
m - 2$ by Proposition 1.2; moreover
$$
\multline
3n + g(E) - m - 2 \le 3(n-1) + g(C) - 1 - m(C)\\
= 3n + g(E) - k - 4 - m(C),
\endmultline \tag1
$$
and the claim follows.

 Let $\gamma$ be the maximal number of nodes of  $E$
approaching $l \backslash
{\bold p}$
as E degenerates to  $C+l$ along various ${\Cal F}(E,C{+}l)$'s.  Set
$$
\Gamma =\sum_{Q} \delta_{Q}(C{+}l) \qquad (Q \in
l\,\cap\,C\backslash
{\bold p}).
$$

{\it Claim\/} 2:\ \  For $\gamma = 0$,  we have  $m = n$,   $\Gamma
= k + 1$,  and the
theorem
follows.  Indeed, by the genus formula \cite{Hi},
$$
\gathered
g(E) - k - 1 = g(C{+}l) = (n-1)(n-2)/2 - \Gamma - \delta_{\bold
p}(C{+}l) -
\sum_{R \notin l} \delta_{R}(C{+}l)\\
g(E) = (n-1)(n-2)/2 - \sum_{Q \rightarrow l \backslash {\bold p}}
\delta_{Q}(E)
-\sum_{Q \rightarrow {\bold p}} \delta_{Q}(E)
-\sum_{Q \not \to l} \delta_{Q}(E)\\
\Gamma-k-1 + \sum_{R \notin l} \delta_{R}(C{+}l)+\delta_{\bold
p}(C{+}l)=
\sum_{Q \rightarrow l \backslash {\bold p}} \delta_{Q}(E)
+\sum_{Q \not \to l} \delta_{Q}(E)
-\sum_{Q \rightarrow {\bold p}} \delta_{Q}(E).
 \endgathered \tag 2
$$
We have
$\sum_{Q \rightarrow {\bold p}} \delta_{Q}(E) \leq  \delta_{\bold
p}(C{+}l)$
and
$\sum_{Q \not\to l} \delta_{Q}(E)\leq
\sum_{R \notin l} \delta_{R}(C{+}l).$
Hence
$$
\sum_{Q \rightarrow l \backslash {\bold p}} \delta_{Q}(E)
\,\ge \,
\Gamma- k - 1\, \ge \, n - 1 - m(C) - k - 1\, \ge \, n - m
$$
(the last inequality follows from Claim 1).  Therefore  $\Gamma \leq
k
+ 1$  and  $\sum_{Q \ne
{\bold p}}(C\cdot l)_{Q}   \leq k + 1.$  Since  $n \geq m,$  we get
$m(C)
= n - k - 2 = m - k - 2,$  and Claim 2 follows.

      We proceed by induction on  $\Delta_{W}.$  Let  $\Cal F (E,C+l)
\subset W$ be a $1$-dimensional family as in the statement of the
theorem.  Since  $C + l$  is a general member of  $W^{\prime},$
one can find a branch
${\Cal W}$  of  $W$  such that  $\Cal F (E,C+l)$, ${\Cal W}^{\prime}
\subset
{\Cal W},$  where  %
${\Cal W}^{\prime} \subset W^{\prime}$  is  a branch through  $C +
l.$
We can now choose a family  $W_{1} \subset U_{m}(n,g+1),$  as
above,
and a branch ${\Cal W}_{1}$ of  $W_{1}$  such that ${\Cal
W}_{1}\supset
{\Cal W}.$
There are two cases to consider.

{\it Case\/} 1:\ \ A general member  $E_{1}^{\prime}$   of some
branch
${\Cal W}_{1}^{\prime}$ of ${\Cal W}_{1} \cap \{a_{0n-m} = 0\}$
containing
${\Cal W}^{\prime}$ does
{\it not\/}  contain
$l.$   Then  $\gamma\ne 0,$  and  $g(E_{1}^{\prime}) = g(E_{1})$
by Proposition 1.2(a).  By a trivial dimension count, some general
member
of           %
${\Cal W}_{1}^{\prime} \cap \{a_{0n-m-1} = 0\}$  is of the form  $C +
l.$
By induction hypothesis, the only singularities of
$C\backslash l$  are nodes and  $m(C) = m + 1 - k_{1} - 2,$  where
$k_{1} = g(E_{1})
 - g(C) = k + 1.$  Therefore  $m(C) = m -k - 2$  and  (1)  is an
equality.

    By induction hypothesis and (2),  $\Gamma - k_{1}-
1 = n - m - 1$  so
$\Gamma - k- 1 = n - m.$  Since
$$
\delta_{\bold p}(C{+}l) = \sum_{Q \rightarrow {\bold p}} \delta_{Q}
(E_{1}^{\prime})
 \leq \sum_{Q \rightarrow {\bold p}} \delta_{Q} (E) \leq
\delta_{\bold p}(C{+}l),
$$
 we get
$$
\sum_{Q \rightarrow {\bold p}} \delta_{Q} (E) = \delta_{\bold
p}(C{+}l),\quad
\sum_{Q \not \to l} \delta_{Q}(E) +
\sum_{Q \rightarrow l \backslash {\bold p}} \delta_{Q}(E) =
\sum_{R \notin l} \delta_{R}(C{+}l) + n - m.
$$

{\it Case\/} 2:\ \ The general members of
{\it all\/}  branches of ${\Cal W}_{1} \cap \{a_{0n-m} = 0\}$  that
contain
${\Cal W}_{1}^{\prime}$ are of the form  $C_1 + l.$   We assume
$\gamma\ne 0.$
There are two possibilities.

 i) $g(C) < g(C_{1}).$  Set $k_{1} =  g(E_{1}) - g(C_{1})$  and  $m(C_{1})
=
(C_{1} \cdot l)_{\bold p}.$  Then  $k_{1} \leq k$  and  $m(C) \geq
m(C_{1}) = m
- k_{1}- 2 \geq m - k - 2.$  Hence  $m(C) =
 m - k - 2$  by Claim 1.  Moreover,  $C$ is a general member of
$U_{m(C)}(n-1,g(C)),$ and
$$
m(C) = m(C_{1}), \  g(C) = g(C_{1}) - 1, \   k = k_{1}, \ \Gamma = n - 1
- m(C),  \
\Gamma - k - 1 = n - m.
$$

We will need the following claim.

{\it Claim\/} 3:\ \  $m(C) =\delta_{\bold p}(C{+}l) - \delta_{\bold
p}(C),$
  and one can assume that at list one of the
following conditions holds:  $m(C) \ne m(C_{1})$,  $\sum_{Q
\rightarrow {\bold p}}
 \delta_{Q}(E) =\delta_{\bold p}(C{+}l)$,  or  $\sum_{Q \rightarrow
{\bold p}} \delta_{Q}(E)
 = 0.$
The first assertion follows at once from the genus formula \cite{Hi}.
Further, if
$\sum_{Q \rightarrow {\bold p}} \delta_{Q}(E)
\ne 0,$  we can assume that
for the general member  $E_{1 }\in {\Cal W}_{1}$,
$$
\sum_{Q \rightarrow {\bold p}} \delta_{Q}(E_{1}) <
\sum_{Q \rightarrow {\bold p}} \delta_{Q}(E)
  \ ( \leq \delta_{\bold p}(C+l) ).
$$
By
induction hypothesis, we have  $\sum_{Q \rightarrow {\bold p}}
\delta_{Q}(E_{1})=
\delta_{\bold p}(C_{1}{+}l)$  hence
$$
\delta_{\bold p}(C_{1}{+}l) <\sum_{Q \rightarrow {\bold p}}
\delta_{Q}(E).
$$
If
 $m(C) = m(C_{1}),$  then  $\delta_{\bold p}(C_{1}{+}l) + 1 \geq
\delta_{\bold p}(C+l) ).$
hence  $\sum_{Q \rightarrow {\bold p}}\delta_{Q}(E) =\delta_{\bold
p}(C{+}l).$
This proves Claim 3.

      If  $\sum_{Q   \rightarrow    {\bold p}}
\delta_{Q}  (E)= 0$
  then  $\sum_{Q \rightarrow {\bold p}} \delta_{Q}(E_1) =
\delta_{\bold p}(C_{1}{+}l) = 0,$  hence  $m(C) = m(C_{1})   = 0$  and
$\delta_{\bold p}
(C{+}l) = 0.$  Therefore  $\sum_{Q \rightarrow {\bold p}}
\delta_{Q}(E) = \delta_{\bold p}
(C{+}l).$  Now, by (2),
$$
\sum_{Q \not \to l} \delta_{Q}(E) +
\sum_{Q \rightarrow l \backslash {\bold p}} \delta_{Q}(E) =
\sum_{R \notin l} \delta_{R}(C{+}l) + n - m.
$$

ii)
$ g(C) = g(C_{1}).$  Then  $k_{1} = k + 1.$  Assume  $m(C) \geq
m(C_{1}) + 1.$  Then
$m(C) > m - k_{1}- 2 = m - k - 3,$  hence  $m(C) =  m
- k - 2$  by Claim 1,  and  $m(C) = m(C_{1}) + 1.$  So  $C$  is moving
in a family of dimension
$$
3(n-1) + g(C_{1}) - 1 - m(C_{1}) - 1 = 3(n-1) + g(C) - 1 - m(C).
$$
Moreover
$\Gamma = n - 1 - m(C)$  and  $\Gamma - k - 1 = n - m.$  If  $
\sum_{Q \rightarrow {\bold p}} \delta_{Q}(E) \ne 0,$  then we can
assume
that  $\sum_{Q \rightarrow {\bold p}} \delta_{Q}(E_1) =
\delta_{\bold p}(C_{1}{+}l)
<\sum_{Q \rightarrow {\bold p}} \delta_{Q}(E).$
It follows from (2) that  $\sum_{Q \rightarrow {\bold p}}
\delta_{Q}(E) = \delta
_{\bold p}(C{+}l),$  and we are done.  But if  $\sum_{Q \rightarrow
{\bold p}} \delta_{Q}(E)
= 0$  then  $\sum_{Q \rightarrow {\bold p}} \delta_{Q}(E_{1}) =
\delta_{\bold p}(C_{1}{+}l) = 0,$  hence  $m(C_{1}) = \delta_{\bold
p}(C{+}l) = 1$
  and exactly one node of  $(C_{1}{+}l) \backslash
{\bold p}$  tends to
${\bold p}$  as  $C_{1}+ l$  tends $ C + l$  (Proposition 1.2(b)).  We
will derive a
contradiction.
 By (2),  $\sum_{R \notin l} \delta_{R}(C{+}l)=\sum_{R \notin l}
\delta_{R}(C_1{+}l).$

      First, we assume that  $C$  and  $C_{1}$
are smooth.  We proceed by induction on
$n - m.$  If  $n = m,$  we degenerate  $C$  into a sufficiently
general nodal curve which is a sum  $\sum
l_{i}$  of  $ n -1$  lines (since  ${\bold p} \in C,$  one of the lines
contains  ${\bold p}$).
Clearly  $\sum
l_{i}$  has  $(n-1)(n-2)/2$  nodes.  Each node of  $\sum l_i$
determines a branch of
${\overline V}_{n,1,0} \subset
{\Bbb P}^{N}$  through
$\sum l_i\, + l.$
We get a degeneration of a reducible nodal curve
$F$  into  $\sum l_i\, + l,$
 where  $F$  has degree  $n$  and contact of order  $n$
with
$l$  at $\bold p;$  moreover,  $F$  contains a line through  $\bold p ,$
which is absurd.  If
$n = m + 1,$  we can split a line
off  $C,$  and apply the preceding argument, etc.

 Next, if  $C \backslash l$  has only unassigned nodes as singularities,
we can
argue as above (note that  $C \backslash l$  may have at most  $n -
2$  such nodes).
But if  $C\backslash l$  contains assigned nodes, we can assume
$\sum_{R \notin l} \delta_{R}(C{+}l)
 \ne \sum_{R \notin l} \delta_{R}(C_1{+}l),$  a contradiction.
This concludes the argument in the case  $m(C) \geq m(C_{1}) + 1.$

Now, we suppose that  $m(C) = m(C_{1})$  for any choice of
${\Cal W}_{1}$  and derive a contradiction.  If  $\sum_{Q \rightarrow
{\bold p}}
\delta_{Q}(E) = \delta_{\bold p}(C{+}l) \ne  0,$  we can assume
$\delta_{\bold p}
(C_1{+}l) \ne
\delta_{\bold p}(C_1{+}l) $   and $\delta_{\bold p}(C_1) \ne
\delta_{\bold p}(C) $  (see Claim 3).  Hence $ g(C) \ne g(C_{1}),$  a
contradiction.

         If  $\sum_{Q \rightarrow {\bold p}} \delta_{Q}(E) = 0,$  then  $
\sum_{Q \rightarrow {\bold p}} \delta_{Q}(E_1) = \delta_
{\bold p}(C_{1}{+}l) = 0$  and  $m(C) = m(C_{1}) = 0$
so  $\delta_{\bold p}(C{+}l) = 0.$  If  $C\backslash l$  contains
assigned points,
we can assume  $\sum_{R \notin l} \delta_{R}(C{+}l)
 \ne \sum_{R \notin l} \delta_{R}(C_1{+}l)$  and derive a
contradiction with a help of (2) and
Claim 1.

  Throughout the rest of the proof, we will exploit that
$
\sum_{Q \rightarrow l \backslash {\bold p}} \delta_{Q}(E)
$
  is {\it too large\/}.  To begin with, we assume, in addition, that
 $C \backslash l$  has no singularities; the case when  $C \backslash
l$  has only
unassigned nodes as singularities is similar.
 By a trivial dimension count, there are two possibilities for
 $C + l$  provided $ \deg(C) \geq 2$:\   $C + l$ has either an ordinary
triple point
or a tacnode.

 First, we consider the case when a general curve  $C_{1}$
  with one node tends to a curve  $C$  with
a node along
$l.$  Then  $d = n - m + 2$  by (2) and Claim 1.  If  $d = 2,$  we
degenerate  $C$
  into a sum $\sum l_i$  of  $n - 1$  general lines with
$l_1 \cap l_2 \cap l \ne \emptyset.$  Then  $\sum l_i$  has  $(n-
1)(n-2)/2 - 1$  nodes outside  %
$l.$  Such a node determines a branch of ${\overline V}_{n,1,0}
\subset
{\Bbb P}^{N}$  through  $\sum l_i \,+l.$  We get a degeneration of
$F$  into
$\sum l_i \,+l,$  where $ F$  is a curve of degree  $n$  having contact
of order  $n$  with  $l$  at  $\bold p$  and  $(n-1)(n-2)/2 - 1 + d$
nodes and no
other singularities.
 Such an  $F$  must contain a line through  $\bold p,$  which is
absurd.

 If  $d = n - m + 2 > 2,$  we degenerate  $C$  into a curve $ C^{\prime}
+ l^{\prime},$
where
$l^{\prime}$  is a general line.
We then split  %
$l^{\prime}$  off the general curve of  ${\overline V}_{n,d,m}$  and
apply the
preceding argument to curves of  ${\overline V}_{n-1,d,m}$
  if  $n - 1 = m,$  etc.

 Next, we consider the case when a smooth curve  $C_{1}$  tends to a
smooth curve
$C$  tangent to
$l$  at one point.  Then  $d = n - m + 1$  by (2) and Claim 1.  If  $d =
1,$  we take
 a general line
$l^{\prime} \subset
{\Bbb P}^{2}$  and a point  $Q' = l^{\prime} \cap l.$  The projection
$\pi_{n,1}$  gives
a natural fibering  $f: \Sigma \rightarrow l^{\prime},$  where
$\Sigma \subset
\Sigma_{n,1,n}$  is an
{\it irreducible\/} subvariety;  $f ^{-1}(Q)$  is a projective space for
every  $Q \in
l^{\prime} \backslash Q^{\prime}.$  By imposing appropriate
conditions
on plane curves, we get a similar fibering  $g: \Sigma ^{\prime}
\rightarrow
l^{\prime}$  whose fibers have dimension  $1.$

       Thus we can assume that  $n = 3,$  and  $g$  is a  $1$-
dimensional
fibering over
$l^{\prime},$ as above.  Then  $g ^{-1}(Q')$  consists of two (genuine)
projective
lines intersecting in one point plus another projective line
of plane curves of degree  $2$  having contact of order  $2$  with
$l$  at  $Q^{\prime}$  and passing through  $2$  general points.  The
latter line will intersect
{\it each\/} of the former lines.  One intersection point corresponds
to a curve
               $l_{1} + l_{2} + l$  with  $Q^{\prime} = l_{1} \cap l_{2} \cap
l,$
and the other one to a curve
$l_{3} + 2l$  with  $Q^{\prime} \notin l_{3}.$  However, there are no
such fiberings over
a projective line, a contradiction.

 If  $d = n - m + 1 >1,$  we split off a line, as in the case $d = n - m +
2 > 2$  above,
and apply the preceding argument provided  $n - 1 = m,$  etc.

It remains to consider the case when  $C$  is a general curve with
unassigned singularities only, one of which is a tacnode or an
ordinary triple point and the rest are nodes  (clearly
$\sum_{Q}\delta_Q(C) \leq n - 2).$
  We will treat the case
when  $C$  has a tacnode and no other singularities; the remaining
cases are similar.  We proceed by induction on $ n - m.$  If  $n = m,$
  we degenerate  C  into a sufficiently general curve  $C^{\prime}$
 which is a sum of  $n - 3$  lines and a quadric tangent to one
of the lines.  As before, we get a degeneration of $ F$  into
$C^{\prime} + l,$  where  $F$
  is a nodal curve of degree  $n$  having contact
of order  $n$  with $l$  at  %
$\bold p;$  moreover,  $F$  contains a line through $\bold p ,$  which
is absurd.  If
$n - 1 = m,$  we can split a line off  $C,$  and apply  the preceding
argument, etc.
This proves the theorem.
\enddemo

\remark {{\rm 2.3.} Remark} $m^{\prime} = \delta_{\bold p}(C{+}l) -
\delta_{\bold p}(C)$
  by the genus formula \cite {Hi}.
\endremark

\remark {{\rm 2.4.} Remark} $m' \leq m - 1$  with equality if and
only if  $E$  contains
a line through $\bold p$ that tends to
$l$  as  $E$  tends to  $C + l.$
\endremark

\head 3. Splitting off a line:  lemmas and a proposition
\endhead

\proclaim{3.1. Lemma}  We keep the notation of Theorem 2.2 and
assume, in addition,
that  $C$  is irreducible and $C \backslash \bold p$    is smooth.
Then  $W_{\text{red}}
^{\prime}\simeq {\Bbb P}^{r}\subset
{\Bbb P}^{N_1}$
where  $N_{1} = (n-1)(n+2)/2$
and  $r = 3(n-1) + g(C) - 1 - m(C).$
Let  $\widetilde U \subset  {\Bbb P}^{r}\   (U \subset {\Bbb P}^{r})$
be the open subset of curves $B$ such that $g(B) = g(C)$
and the corresponding multiplicities are equal:  $b_{j}(B) = b_{j}(C) =
b_{j}(W^{\prime})$
(in addition,  $B$ intersects  $l \backslash \bold p$
transversely).  We set $\Gamma = C  \cap (l\backslash \bold p)$
and assume  $\text{\rm length}(\Gamma) = \deg(C) - m(C) \geq 2.$
Consider the incidence correspondence  $I \subset (l\backslash\bold
p)\times U$
and the projection $I \rightarrow U.$
 Then the image of the monodromy map
$$
\mu : \pi_1(U,C) \rightarrow \text{\rm Aut} (\Gamma)
$$
is the full symmetric group.
\endproclaim

\demo {Proof}  Since the multiplicities $ b_{j}(W^{\prime})$'s
determine  $g(C)$,
$W_{\text{red}}^{\prime}$  is a linear system by
Proposition 1.2 and its proof.  We set
$$
I(2) =\{(p _{1},  p _{2}, D) \bigm |   p _{1},  p _{2} \in
(l\backslash\bold p) \cap D,\,   p_1\ne p_{2}\} \subset l \times l
\times
{\Bbb P}^{r}
$$
where
 $g(D) = g(C)$  and  $b_{j}(D) = b_{j}(C)$  for all  $j.$  The
image of  $\mu $  is the full symmetric group, provided it is twice
transitive and contains a simple transposition.  As in \cite{ACGH, pp.
111 -- 112}, it will
suffice to verify that the fibers of the projection  $I(2)  \rightarrow
(l\times l
)\backslash \Delta$  are projective spaces and  $\widetilde U$
 contains curves simply tangent to
$l\backslash \bold p$ at one point.  Both properties follow at once
from
Proposition 1.2 and its proof.
\enddemo

\proclaim{3.2. Lemma}   Let  $W^{\prime}$  be a component of
${\overline V}_{n,d,m}
\cap\{a_{0n-m} = 0\}$  with a general member $D = C + l$  such that
$\delta_
{\bold p}(C) = 0.$   Then ${\overline V}_{n,d,m}$ contains a branch
through $(D, b)$
 with a cycle $b$ containing all nodes of $C.$
\endproclaim

\demo {Proof}  First, one can find a connected cycle on  $D$  of the
form
$$
 b = \delta_{\bold p}(C)\cdot \bold p + b_{1} + b_{2},$$
 where
$b_{1}$  is a sum of  $n - m$  nodes of  $D$  along
$l\backslash \bold p,$  and  $b_{2}$  is a sum  of the nodes of  $C.$
Indeed,
consider {\it any} branch of  ${\overline V}_{n,d,m}.$  We get a cycle
of assigned
singularities on  $D:  \delta_
{\bold p}(C)\cdot \bold p + c_{1} + c_{2},$  where  $c_{1}$  is a sum of
nodes of  $D$  along      %
$l\backslash \bold p.$
By Theorem 2.2,  $C$  has unassigned nodes if and only
if  $c_{1} > n - m.$

 If  $c_{1} > n - m,$  we proceed as follows.  Let $C^{\prime}$
  be a component of  $C$  with an unassigned point  $Q_{1}$  on $l.$
If all nodes of  $D$
  along  $C^{\prime}\backslash l$  are assigned, we consider  $D -
C^{\prime}$
  in place of  $D.$     Otherwise, let  $Q_{2}$  be an unassigned node of
$D$  along
$C^{\prime}\backslash l,$  and  $Q_{2}$  also belongs to a component
$C^{\prime\prime}.$
  If  $C^{\prime\prime}$  has an assigned point  $Q_{3}$  along
$l\backslash \bold p,$
we can interchange  $Q_{2}$  with  $Q_{3}.$  If  $C^{\prime\prime}$
has an
unassigned point  along
$l\backslash \bold p,$  we can interchange  $Q_{2}$  with any
assigned node of
$D \backslash \bold p.$  Now assume  $C^{\prime\prime}$  is a line
through  $\bold p.$
If  $C^{\prime\prime}\backslash Q_{2}$  has no unassigned points,
we consider
$D - C^{\prime\prime}$  in place of  $D.$  Finally, suppose
$C^{\prime\prime}\backslash
Q_{2}$  has an unassigned node  $Q_{4},$  and  $Q_{4}$  also belongs
to a component  $F.$
  Then we
interchange  $Q_{4}$  either with $ F\cap l$  (if  $ F\cap l$ is
assigned) or with any
assigned node of  $D.$

Now, let  ${\Cal W}$  be a branch of  ${\overline V}_{n,d,0}$  through
$(D, b).$
We claim that ${\Cal W}\cap\{a_{0n} = \dots = a_{0n-m+1} = 0\}$
contains the required
 branch of  ${\overline V}_{n,d,m}.$  Indeed, if a general member of
${\Cal W}\cap\{a_{0n} = \dots = a_{0n-m+r} = 0\},  \ r \geq 2,$
has the form  $F + l,$  then  $F$  and  $C$  have the same number of
nodes,  a
contradiction.
\enddemo

\proclaim{3.3. Lemma}   For $0 \leq g \leq (n-1)(n-2)/2 $  and  $0
\leq m \leq n,$
 there exists a unique component of $U_m(n,g),$  denoted by
${\overline V}_{n,d,m} \
  (d = (n-1)(n-2)/2 - g),$  whose general members are smooth at $
\bold p$   and irreducible.
\endproclaim

\demo{Proof}  The existence is a classical fact and the uniqueness
is proved, for instance, in  \cite{R, Irreducibility Theorem (bis)}.
 We will present a proof of the existence of  ${\overline V}_{n,d,m}.$

 For  $m = 0$  or  $1,$  the existence is well known and follows
from the deformation theory.  We assume  $m \geq 2.$  Let  $D + l$
be a curve, where  $D$
  is a sum of  $n - 1$  general lines.
 We choose  $d - n + m$  nodes  of  $D$  and  $n - m$  nodes of    $D +
l$  along
$l$  such that  $D +l$  remains connected after blowing up these
nodes.  Regarding
these  $d$  nodes as assigned, we get a family of irreducible curves
with  $d$  nodes and no other singularities in a neighborhood of $D
+l.$  Intersecting this
family with  $m$  hyperplanes,  $a_{0n}
= \cdots = a_{0n+1-m} = 0,$  we derive the existence of
a required component by Theorem 2.2.
\enddemo

\proclaim{3.4. Lemma} Let $W$  be a component of  $  U_m(n,g),$
and $E$  a general
member of $W.$
Let $W^{\prime}$ be a component of $W \cap \{a_{0n-m} = 0\}$
and $D$
a general member of $W^{\prime}.$
If   $l \nsubseteq
 D,$    then $g(D) = g,$ and  $\sum_{Q \rightarrow {\bold p}}
\delta_{Q}(E)=
\delta_{\bold p}(D)$  as  $E$ tends to $D.$
If, in addition,  $\delta_{\bold p}(D) = 0$ then  $W$  is smooth at $D.$
\endproclaim

\demo{Proof}  The first assertion follows from Proposition 1.2
and the genus formula \cite{Hi}.  Further, if  $E$  is smooth, then
 $W^{\prime} = W_{\text{red}}^{\prime}$  and  $W$  is smooth at
$D.$  But if  $E$
  has nodes, we get
$
{\Cal O}_{D,W^{\prime}} = (
{\Cal O}_{D,W^{\prime}}
)_{\text{red}}$  by a simple induction on the
number of nodes of  $E,$  hence  $W$  is smooth at  $D$  by
Proposition 1.2.
\enddemo

\proclaim{3.5. Lemma}  With the notation of Lemma 3.2,
$$
\dim{\bold T}_{D}({\overline V}_{n,d,m})  \geq  \dim {\overline
V}_{n,d,m} +n-m{.}
$$
\endproclaim

\demo{Proof}  For simplicity, we assume  $C$  is irreducible and
$n\ne m.$  Each branch  %
${\Cal W}$ of  ${\overline V}_{n,d,m}$  through  $D$  determines the
corresponding
cycle of assigned singularities of  $D,$  denoted by  $c({\Cal W} ).$
By Lemma 3.2,
${\overline V}_{n,d,m}$  contains a branch,
denoted by  ${\Cal W}_{1},$ through  $(D, b)$  with  $c({\Cal W}_{1})$
containing all nodes
of  $C.$  Then  $c(
{\Cal W}_{1})$  contains exactly  $n - m$  nodes of  $D$
along  $l\backslash \bold p.$

 Now, we exhibit  $n - m + 1$  branches
${\Cal W}_{1},\dots, {\Cal W}_{i}$  such that the span of the tangent
spaces to
${\Cal W}_{1}, \dots ,{\Cal W}_{i-1}$  does not contain the tangent
space to
${\Cal W}_i$, $i = 2, \dots,   n - m + 1.$  We take  $n - m + 1$
 nodes of  $D\backslash \bold p,$ containing  $n - m$  assigned (with
respect to  %
${\Cal W}_{1})$  nodes.  Such nodes exist by the Principle
of Degenerations because  $\sum_{Q \rightarrow {\bold p}}
\delta_Q(E) = \delta_
{\bold p}(D).$  By standard deformation theory, we obtain $n - m +
1$  branches  $\{
{\Cal W}_{i}\}$  of  ${\overline V}_{n,d,m}$  with
$$
 \sup c({\Cal W}_{1}) \cap \dots \cap \sup c(
{\Cal W}_{i-1}) \nsubseteq \sup c({\Cal W}_{i}) \quad   (i = 2, \ldots ,
n-m+1)
$$
 and
the desired tangent spaces. (Note: if  $C$  is reducible, then
{\it each\/} component of  $C$  contains unassigned nodes of  $D$
along
$l.$)
\enddemo

\remark{{\rm 3.6.} Remark}  With the notation of Lemma 3.2, we
assume
that  $\delta_{\bold p}(D) = 0.$  Let  %
${\Cal W}$   be a branch of  ${\overline V}_{n,d,m}$  through  D  with
 $n - m$  assigned nodes along  %
$l$  and  $d-n+m$  assigned nodes away from
$l.$  Then
${\bold T}_{D}(%
{\Cal W} ) \subset {\Bbb P}^N$ consists of all curves of degree   $n$
having
contact of order at least  $m$  with $l$  at  %
$\bold p$  and passing through all the assigned nodes of  $D$ (cf.
\cite{H,~Sect.~2}, \cite
{Z});  moreover,  $\{a_{0n-m} = 0\}$  intersects
${\Cal W}$   transversely.
\endremark

\subhead {3.7} \endsubhead  By abuse of notation, we will denote by
$\pi_{N}$  and
  $\pi_{n,d}$
the restrictions of the corresponding projections to  $\Sigma_{n,d,m}$
and its tangent
spaces.

 The following lemma is a special case of Proposition 3.9 below.%

\proclaim{3.8. Lemma}  For $2 \leq m \leq n$ and
$n - m + 1 \leq d \leq (n - 1)(n - 2)/2,$ the scheme
 ${\overline V}_{n,d,m} \cap \{a_{0n-m} = 0\}$
contains a component  $W^{\prime}$ whose general members have
the form
$D = C + l,$ where  $m(C) = 1$ and $C$ has $d - 1 - n + m$ nodes.
Moreover, there exists
$
(D,\bold p + b) \in \Sigma_{n,d,m}$ such that $b$ contains all the
nodes of $C$ and
$$
 \dim
{\bold T}_{(D,\bold p + b)}(\Sigma_{n,d,m})  =  \dim W^{\prime} + 2.
$$
\endproclaim

\demo{Proof}  {\it Case}: \ $d = n - m + 1.$  First, we assume that
$n = m$  hence  $d = 1.$  By a trivial dimension count, a general
member of a component of  ${\overline V}_{n,1,n} \cap \{a_{00} =
0\}$
is a curve of the form  $F +l,$  where  $F$  may be either a smooth
curve through
$\bold p,$
  or a curve with one node.  By Theorem 2.2,  ${\overline V}_{n,1,n}
\cap \{a_{0n-m} = 0\}$
 contains the component  $K$  whose general members have the
form
$F + l$  with
$F \in U_{0}(n-1,(n-2)(n-3)/2 - 1).$  Indeed,  ${\overline V}_{n,1,n}
\cap \{a_{0n-m} = 0\}$
 contains a member that is a sum  $\sum l_i+l$  of  $n$  distinct lines
meeting in a point
of  $l\backslash \bold p$ \cite{Z, Lemma 2}.  In fact, for  $0 \leq m
\leq n,$  $1
\leq d \leq (n - 1)(n - 2)/2,$  and an arbitrary point of  %
$l$,  ${\overline V}_{n,d,m}$ contains all curves of the form
$\sum l_i+ l$ that are sums
of  $n$  distinct lines meeting in that point (see (3.3)).

Now, we consider  $\Sigma_{n,1,n-1} $ and its point  $\alpha = (\sum
l_i + l, \bold p),$  where  %
$l_{i}$  are general lines meeting in  %
$\bold p \ (1 \leq i \leq n - 1).$  Set  $H_{0}= \pi_N^{-1}(\{a_{01}=
0\}).$  In a small
neighborhood
of  $\alpha$  in  $\Sigma_{n,1,n-1} $, $\Sigma_{n,1,n-1} \cap H_{0}$
is an analytic subset,
denoted
by  ${\Cal W},$  connected in codimension $1.$  This follows at once
from Proposition 1.3 and Grothendieck's connectedness theorem
\cite{G, Exp XII, Theorem 2.1} (compare  \cite{Tr, Lemma 1}).  Let
$U$  be the closure in  $
{\Bbb P}^{N}$  of  $U_{n}(n,(n-1)(n-2)/2 - 1) \backslash {\overline
V}_{n,1,n}.$  Clearly  $U$
  is a linear system.
 We can describe the branches of  ${\Cal W}_{\text{red}}:$

 a) a (unique) branch
${\Cal A}$ of  $\Sigma_{n,1,n} $ (see Lemma 3.3 and Proposition 1.3);

 b) the branches of  $\pi_N^{-1}(U)$  inside  ${\Cal W}_{\text{red}}$;
and

 c) a (unique) branch whose general point is  $(E+l, c),$  where  $E$
is a general curve of
degree  $n - 1$
 through  $c$\ ($c\ne \bold p$).

Since  $U$  is a hyperplane in  ${\overline V}_{n,0,0}$,  $H = U \cap
{\overline V}_{n,1,n}$
is a hypersurface in the projective
space  $U.$  We will show that  $H_{\text{red} }$ is the required
$W_{\text{red}}^{\prime}.$
 By a trivial dimension count, it will suffice to show that
$H_{\text{red}}$  contains a component that lies in $\{a_{00} = 0\}.$
Suppose to the
contrary that all general points of
$ H \cap\{a_{00 }= 0\}$  are curves of the form $F^{\prime} + l$,
where  $F^{\prime}$  is a
general curve of degree  $n - 1$  through
$\bold p$  with one node.  Let  $K^{\prime}$ denote the
corresponding component
of  $H \cap\{a_{00 }= 0\}.$  The key point is that
$K_{\text{red}}^{\prime}$
 is {\it not\/}  a linear system.  We get
$$
 \dim H = \dim K^{\prime} + 1 = n(n + 3)/2 - n - 2 = n(n - 1)/2 + n -
2.
$$

First, we will show that each component of  $H_{\text{red}}$  is a
hyperplane in  $U.$
Let  $E^{\prime} = l_{1}^{\prime}+\cdots+ l_{n-1}^{\prime}+ l \in
K^{\prime}$  be a nodal
curve, where  $l_{1}^{\prime}, \dots , l_{n-1}^{\prime}$  are general
lines and
$\bold p \in  l_1^{\prime}.$
  We choose one general point on each
$l_i^{\prime}  \ (1 \leq i \leq n - 1),$  and consider the set
$M$  consisting of these points and the nodes of
$E^{\prime}\backslash \bold p.$
 Thus  $M$  has  $\dim H$  points.  To each point, there
corresponds a hyperplane in  $
{\Bbb P}^{N}$  of curves passing through that point.  Clearly
those hyperplanes will intersect  $K^{\prime}$  in a %
{\it reduced\/} scheme consisting of one point, namely  $E^{\prime}.$
On the other hand, if a component  $H^{\prime} \subset
H_{\text{red}}$  is %
{\it not\/} a hyperplane, then  $H^{\prime} \cap\{a_{00}= 0\}$
intersects
those  $\dim H$  hyperplanes in a subscheme containing more than
one point, because  $H^{\prime} \cap\{a_{00}= 0\}$  belongs to
the hyperplane corresponding to the point $l_{n-1}^{\prime}\cap l.$

 Now take any component of  $H_{\text{red}}$,  say  $H^{\prime}
\subset U.$  Clearly
 $\dim (H^{\prime} \cap\{a_{00}= 0\}) \geq \dim (H^{\prime} \cap
K^{\prime})$
 and we derive a contradiction.  We have also proved that
$$
H_{\text{red}} = U \cap {\overline V}_{n,1,n} \cap\{a_{00 }= 0\} =
U \cap\{a_{00 }= 0\} = W_{\text{red}}^{\prime}.
$$

We set  $F = \sum l_i,$  where  $l_{i}$  are general lines meeting in
$\bold p
 \ (1 \leq i \leq n - 1),$ and take an additional general line
$l^{\prime}$  through  $\bold p.$  To obtain two lines in
$\Sigma_{n,1,n} $
through  $\alpha,$
  we move  $F$  along the
two lines, namely  $l$  and  $l^{\prime}.$  Thus
$$
\pi _{n,1}: \
{\bold T}_{\alpha}(\Sigma_{n,1,n} ) \rightarrow {\bold T}_{\bold
p}({\Bbb P}^{2})
$$
 is
surjective hence  $\dim {\bold T}_{\alpha}(\Sigma_{n,1,n} ) \geq
\dim W^{\prime} + 2.$
  Since
$$
{\overline V}_{n,1,n}  \cap\{a_{00} = a_{1n-1}= 0\} =
W_{\text{red}}^{\prime},
$$
$\{a_{00} = 0\}$  is transversal to  ${\overline V}_{n,1,n}$  at  $F + l,$
and  $\{a_{1n-1}= 0\}$
is transversal to    ${\overline V}_{n,1,n}  \cap\{a_{00} = 0\}$
at  $F + l.$  Therefore  $a_{00  }$ and $a_{1n-1}$  produce linearly
independent elements in
${\bold T}_{(D,\bold p)}({\overline V}_{n,d,m})$  and  $\dim {\bold
T}_{(D,\bold p)}
(\Sigma_{n,1,n} ) = \dim W^{\prime} + 2.$

Now assume  $n\ne m.$  We proceed by induction on  $d.$  Consider
 ${\overline V}_{n,d-1,m} \cap\{a_{0n-m} = 0\}.$
 By Lemma 3.3, it contains  ${\overline V}_{n,d-1,m+1}.$  By
induction hypothesis,
${\overline V}_{n,d-1,m+1} \cap\{a_{0n-m-1} = 0\}$
          contains
a component  $W^{\prime}$  that is also the required component of
${\overline V}_{n,d,m}
\cap\{a_{0n-m} = 0\}$  (Theorem 2.2).  Indeed, intersect
an appropriate branch of  ${\overline V}_{n,d-1,m}$  with a branch
of a hypersurface
 ${\overline V}_{n,1,0}\subset
{\Bbb P}^{N}$  corresponding to an unassigned node along
$l$  of a general member of  $W^{\prime}.$  The remaining assertions
of the lemma also follow at once (see Lemma 3.2).

{\it General case}.  We consider  ${\overline V}_{n,d^{\prime},m}$
with  $d^{\prime}
= n - m + 1.$  One can find a component of  ${\overline
V}_{n,d^{\prime},m}\cap \{a
_{0n-m} = 0\}$  with a general member  $D^{\prime} = C^{\prime} +
l$
that satisfy the lemma.  We degenerate  $C^{\prime}$  into a
nodal curve  $C$  with  $d - d^{\prime}$  nodes and  $m(C^{\prime})
= m(C)$  (Lemma
3.3),  and obtain a curve  $D = C + l$  that is a general member of the
required component
of  ${\overline V}_{n,d,m} \cap\{a_{0n-m} = 0\}.$  The remaining
assertions of the lemma also follow at once (see Lemma 3.2).
\enddemo

\proclaim{3.9. Proposition}   For $m\leq n,$  let  $W^{\prime}$

be a component of ${\overline V}_{n,d,m} \cap\{a_{0n-m}
= 0\}$  with a general member   $D = C + l = C_{1} + \cdots + C_{q} +
l,$  where each
$C_{r}$  is irreducible  $(1 \leq r \leq q).$
Let $  a = e\bold p+ b\in \text{\rm Sym}^{d}({\Bbb P}^{2}),\
\bold p\notin \sup(b),$
be a cycle of assigned singularities of  $D$  (with respect to a branch
of
${\overline V}_{n,d,m}$) such that $b$ contains all nodes of  $D
\backslash l.$  Then

{\rm i)}\ \ $\dim {\bold T}_{(D,a)}(\Sigma_{n,d,m} ) = \dim
W^{\prime} + e + 1;$

 {\rm ii)}\  $\dim {\bold T}_{D}({\overline V}_{n,d,m}) \geq \dim
{\overline V}_{n,d,m} + e
+ n - m;$ and

{\rm iii)} $\delta_{\bold p}(C) = 0,$ that is,  $C$ is smooth at $\bold
p$.

Moreover, for every integer  $e$,  $0 \leq e \leq \min \{d
- n + m,  m - 2\}, $ there exist components  $W^{\prime}$   and
points  $(D, a) \in \Sigma
_{n,d,m}$  as above.
\endproclaim

\demo{Proof}  For  $e = 0,$  the proposition  follows at once
from Theorem 2.2 and Lemma 3.3.  In fact, there is then only
one (smooth) branch of  $\Sigma_{n,d,m}$ through  $(D, a).$  We
assume that  $e \geq 1$  and proceed by induction on  $e.$  We have
two natural inclusions:
$$
{\Bbb P}^{e} = \text{Sym}^{e}(l) \subset \text{Sym}^{e}({\Bbb P}^{2}),
\qquad
{\bold T}_{e{\bold p}}(\text{Sym}^{e}({\Bbb P}^{2})) \subset {\bold
T}_{a}
(\text{Sym}^{d}({\Bbb P}^{2})),
$$
 and
a natural map  ${\bold T}_{(D,a)}(\Sigma_{n,d,m}) \rightarrow
{\bold T}_{e{\bold p}}(\text{Sym}^{e}({\Bbb P}^{2})).$

{\it Case\/}:  $\sup(b) \subset l.$  We observe that
${\bold T}_{(D,a)}(\Sigma_{n,d,m}) \subset {\bold
T}_{(D,a)}(\Sigma_{n,d,0})$,
$\Sigma_{n,d,0 }$ is smooth at  $(D, a),$  and  $\pi_{N  }$ induces a
one-to-one map of
a small neighborhood
of  $(D, a)$  in  $\Sigma_{n,d,0}$ onto its image in  ${\overline
V}_{n,d,0}$
(\cite{AC,~Sect.~4},    \cite{DH2,~Sect.~4}, \cite{Ta}).  So $\pi_N$
induces an embedding
${\bold T}_{(D,a)}(\Sigma_{n,d,m}) \subset {\bold T}_{D}({\overline
V}_{n,d,m}).$
  By Lemma 3.5,  ${\overline V}_{n,d,m}$
 may have several branches through  $D$  (at least if  $\delta_{\bold
p}(C) = 0$  and
$n > m),$  and we will see, among other things, that
$$
 \dim {\bold T}_{D}({\overline V}_{n,d,m})  \geq  \dim {\bold
T}_{(D,a)}(\Sigma_{n,d,m})
 + n - m.
$$
By Theorem 2.2 and Remark 2.3,
$$
m(C) = m - (g - g(C)) - 2  \leq  e = \delta_{\bold p}(C{+}l)
$$
 with
equality only if  $\delta_{\bold p}(C) = 0.$  By induction hypothesis,
we get  $d - 1
\leq n - 2.$  We will also show that
$$
 W^{\prime}
\cap\{a_{0n-m-1} = \cdots =a_{00 }= a_{1n-1}= \cdots = a_{1n-e}=
0\} =
W_{\text{red}}^{\prime}.
$$

Consider a general point of  $\Sigma_{n,d,m  }$ that is very near
 $(D, a).$  Regarding one of its nodes approaching $\bold p$ as
virtually non-existent,
we get a component  $W^{\prime\prime}$
 of  ${\overline V}_{n,d-1,m} \cap\{a_{0n-m} = 0\}$
 such that  $W^{\prime} \subset W^{\prime\prime}$  and a general
member of
$W^{\prime\prime}$ contains
$l$  (Theorem 2.2).  By induction hypothesis
$$
\split \dim W^{\prime\prime} \cap\{a_{0n-m-1} = \cdots =a_{00 }=
a_{1n-1}=&
\cdots = a_{1n-e+1}= 0\}\\
&  = \dim W^{\prime\prime} \,( = \dim W^{\prime} + 1).
\endsplit
$$
 Since
$$
 \split W^{\prime}\cap\{a_{0n-m-1}=\cdots = a_{1n-e }= 0\}
\subset W^{\prime\prime}\cap
\{a&_{0n-m-1 }=\cdots = a_{1n-e}= 0\}\\
& = W_{\text{red}}^{\prime\prime}\cap
\{a_{1n-e }= 0\},
\endsplit
$$
 we get
$$
 \dim W^{\prime}\cap\{a_{0n-m-1 }= \cdots = a_{1n-e }= 0\}  \leq
\dim W^{\prime}.
$$
 This
inequality is an equality only if  $W^{\prime} \cap \{a_{0n-m-1
}=\cdots =
 a_{1n-e}= 0\} = W_{\text{red}}^{\prime}.$

Now we assume, in addition, that  $e \geq 2$  and  $n = m,$  and
generalize  Lemma 3.8 where the case  $e = 1$  was discussed.
 Consider  $\Sigma_{n,e,n-1}$ and its point  $\alpha_0 = (nl, e\bold
p).$  Set  $H_0 =
\pi_N^{-1}(\{a_{01}= 0\}).$  Let  $ U$  be the closure in
${\Bbb P}^{N}$  of
$$
U_{n}(n, (n-1)(n-2)/2 - e) \backslash {\overline V}_{n,e,n}.
$$
In a small
{\it neighborhood\/} of  $\alpha_0$  in  $\Sigma_{n,e,n-1} $,
$\Sigma_{n,e,n-1}
\cap H_0$  is an analytic subset, denoted
by $\Cal W ,$  connected in codimension  $1.$  We can describe the
branches of  $\Cal W_{\text{red}}:$

 a) a (unique) branch ${\Cal A}$ of  $\Sigma_{n,e,n}$ (see Lemma 3.3
and Proposition
1.3);

 b) the branches of  $\pi_N^{-1}(U)$  inside {$\Cal W_{\text{red}}$;
and

 c) for each  $t$,  $0 \leq t \leq e-1,$  the branches, denoted
by ${\Cal A}_{t}$,  ${\Cal A}_t^{\prime}, \dots ,$  whose general
points have
the form  $(E+l, t\bold p +c),$  where  $\bold p\notin \sup(c)$,
$\delta_{\bold p}(E) = 0$,
and E  is a general member of a component of  $U_{t}(n-1, (n-2)(n-
3)/2 - e + 1 + t).$

In (c), we used the induction hypothesis;  note that  $(E+l, t\bold p
+c)$  have at least
{\it one assigned\/} node along $l\backslash \bold p$  (Theorem
2.2).

 To prove (iii), we consider the branch ${\Cal W }^{\prime\prime}$
of  $W^{\prime\prime}$
  whose general member  has the form  $D^{\prime\prime} =
C^{\prime\prime} + l$  with
$m(C^{\prime\prime}) = e - 1.$  Let
${\Cal V}$  be a branch of  $\pi_{N}({\Cal A}) \cap \{a_{00 }= 0\}$
whose general member
is the curve  $D = C + l.$  We may assume that  $(D^{\prime\prime},
(e-1)\bold p)$  tends
to  $(D, (e-1)\bold p)$  (note that  $\delta_{\bold p}(D) = e$  while
$\delta_{\bold p}
(D^{\prime\prime}) = e - 1).$

 We claim that ${\Cal V}_{\text{red}} \subset \pi_{N}({\Cal A}_{e-
1})$  for a suitable
${\Cal A}_{e-1},$  and $(E+l, (e-1){\bold p}+c)$  tends to  $(D, a).$
Then  $c$
will approach $\bold p,$  and we get  $m(C) = e$  hence  $C$  is a
smooth curve.
 To prove the claim, we observe that
$$
\pi_{N}(\Sigma_{n,e,n-1}) \subset \pi_{N}(\Sigma_{n,e-1,n-1}),
\qquad
{\Cal W }^{\prime\prime} \subset \pi_{N}(\Sigma_{n,e-1,n-1}) \cap
\pi_{N}(H_0).
$$
 As
before, we can describe the branches of $\pi_{N}(\Sigma_{n,e-1,n-1})
\cap \pi_{N}(H_0).$
  By induction hypothesis, these
branches have general members of the form $ E^{\prime} + l$  with
$m(E') \leq e - 2,$  provided those members
contain  $l.$  By a trivial dimension count, only the branch with
 $m(E^{\prime}) = e - 2$  can contain ${\Cal
W}_{\text{red}}^{\prime\prime}$
  and  $\pi_{N}
({\Cal A}_{e-1})_{\text{red}}.$  Therefore  ${\Cal
W}_{\text{red}}^{\prime\prime}
= {\pi_{N}
({\Cal A}_{e-1})}_{\text{red}}.$  So $
{\Cal V}_{\text{red}} \subset \pi_{N}({\Cal A}_{e-1}).$  Note that
$W_{\text{red}}^{\prime}$
is a linear system.

Next, we will prove (i) and (ii) in case  $n = m$  (that is,  $a = e\bold
p$).  Let
$$
 \alpha = (l_1+ \cdots + l_{n-2}+ 2l, e\bold p) \in \Sigma_{n,e,n} ,
$$
 where $l_{1},\dots  , l_{n-2}$  are general lines meeting in  $\bold
p.$  Moving
$l_{1}$  along  $l,$  we get a point
$$
(l_1^{\prime}+l_2+ \cdots + l_{n-2}+ 2l, \Sigma p_i) \in \Sigma_{n,e,n}
,
\quad p_1=l_1^{\prime} \cap l, \quad p_i=l_1^{\prime} \cap l_i\ (2
\leq i \leq e).
$$
We then take arbitrary  $e$  points on $l$:  $q_{1},\dots ,q_{e}.$
Consider  $e$  general lines
$l_i^{\prime}$   with  $q_{i }=l_i^{\prime} \cap l \  (1 \leq i \leq e).$
Then
$$
(l_1^{\prime}+\cdots+l_e^{\prime}+l_{e+1}+ \cdots + l_{n-2}+ 2l,
\Sigma q_i)
 \in \Sigma_{n,e,n};
$$
see the proof of Lemma 3.3 (existence).  This shows that the image
of the natural map
$$
{\bold T}_{\alpha}(\Sigma_{n,e,n} )  \rightarrow
{\bold T}_{e\bold p}(\text{Sym}^{e}({\Bbb P}^{2}))
$$
contains  ${\bold T}_{e\bold p}(\text{Sym}^{e}(l))$
and has dimension at least  $e + 1.$  It follows that
$$
\dim {\bold T}_{\alpha}(\Sigma_{n,e,n}) \geq \dim W^{\prime} + e +
1.
$$
Since
$$
{\overline V}_{n,d,n}\cap \{a_{00}= a_{1n-1} = \cdots =a_{1n-e}=0 \}
 = W_{\text{red}}^{\prime},
$$
 we
get  $\dim {\bold T}_{(D,e\bold p)}(\Sigma_{n,e,n} ) \geq \dim
W^{\prime} + e + 1,$
  as in the
proof of Lemma 3.8.  This inequality, however, must be an equality.

To prove (i) - (iii) in case  $n\ne m,$  we proceed by induction
on the number of assigned nodes along  $l\backslash \bold p,$  as in
the proof of
Lemma 3.8 (case: $n\ne m$).  Note
that if  $\sum  a_{jk}X^{j}Y^{k}Z^{n-j-k} = 0$  is an equation of a
general
member of  ${\overline V}_{n,d,m},$  then  $a_{0n-m }$, \dots ,
$a_{00 }$, $a_{1n-1}$, \dots , $a_{1n-e}$  produce  $n - m + e + 1$
linearly
independent elements in  $
{\bold T}_{D}({\overline V}_{n,d,m}).$ %

To establish the existence of  $W^{\prime},$  we generalize the
corresponding
argument from Lemma 3.8.  We assume  $n = m,$  because the case
 $n\ne m$  will follow by induction on  $d,$  as in Lemma 3.8.
 Consider the branches of  $\Sigma_{n,e,n} \cap \pi_N^{-
1}(\{a_{00}=0\})$
  through the point $\alpha_{0}.$  It is known what kind of branches
one may expect
(Theorem 2.2).  By induction hypothesis, for each  $t$,  $0 \leq t \leq
e-1,$  we get branches whose general points have the form
$(E+l,t\bold p+c),$
  where  $\bold p\notin \sup(c)$,  $\delta_{\bold p}(E) = 0,$  and $E$
is a general
member of a component of  $U_{t}(n-1, (n-2)(n-3)/2 - e + t).$  We
claim that for
$t = e,$  we get similar branches.

 Let  $\gamma (e,t)$,  $0 \leq t \leq e$,  denote the number of
conditions
imposed on curves of {\it sufficiently large\/} degree to have contact
of order at
least  $t$  with  $l$  at  $\bold p$  and a singularity at  $\bold p$
with
$\delta_{\bold p} \geq e - t$  \cite{S, Lemma 6}.  Here the general
members
of the ambient variety are curves, say of degree  $r,$  having
smooth contact of order  $t$   with  $l$  at  $\bold p$ and  $e - t$
nodes and no
other singularities;  we also consider the corresponding general
points
$(F, t\bold p+c) \in {\Bbb P}^{R }\times \text{Sym}^{e}({\Bbb P}^{2}),$
where  $R = r(r + 3)/2$  and $\bold p\notin \sup(c).$

 First, we assume that  $n - 1$  is sufficiently large.  Thus the
dimension of the subfamily of the corresponding ambient variety,
whose general members have degree  $n - 1,$  contact of order at
least  $t$  with $l$  at  $\bold p,$  and a singularity at $\bold p$
with  $\delta_{\bold p}
 \geq e - t,$  equals
$$
(n - 1)(n + 2)/2 - e - \gamma (e,t) = n(n + 3)/2 - e - n - \gamma (e,t)
- 1 \quad
  (0 \leq t \leq e).
$$

For  $0 \leq t \leq e - 1,$  let  $\Sigma_t \subset
{\Bbb P}^{N} \times  \text{Sym}^{e}({\Bbb P}^{2})$  be the closed
subvariety of
the  corresponding
ambient variety, whose general members have the form  $(F, e\bold
p),$  where  $F$
  has degree  $n,$  contact of order at least $t$  with $l$  at  $\bold
p,$
and a singularity at  $\bold p$  with  $\delta_{\bold p} \geq e - t.$
  In particular  $F$  has a singularity at  $\bold p.$
Utilizing the proof of the existence in Lemma 3.8, we get
$$
\dim (\Sigma_{n,e,n}\cap \pi_N^{-1}(\{a_{00}=0\}) \cap \Sigma_t)
 = \dim (\Sigma_{n,e,n} \cap \Sigma_t)
  \geq  n(n + 3)/2 - e - n - \gamma (e,t),
$$
 where
 $\dim (\cdot)$  means, as usual, the dimension of the components of
maximal dimension.  Comparing this estimate with the previous
one for  $t \leq e - 1,$  we obtain a required  $W^{\prime}.$

 Finally, if  $n - 1$  is not sufficiently large, we take a sufficiently
large integer  $r$  and apply the preceding argument to
$\Sigma_{r+1,e,r+1} $
in place of  $\Sigma_{n,e,n}.$  We then split   $r - n + 1$ {\it
general\/}  lines off.

{\it General case.}  As we have already seen before, the nodes
of  $D\backslash l$ do not play an essential role.  The existence is
established
by induction with the help of Lemma 3.3 (uniqueness); see the
corresponding
argument in Lemma 3.8.  To prove (iii), we do not
need a generalization of Proposition 1.3.  Assume  $D\backslash l$
has
$v$  nodes.  Regarding those nodes as virtually
non-existent, we apply the above discussion to  ${\overline
V}_{n,d',m} \
(d^{\prime} = d - v)$  and the original curve  $D.$  We now consider
$\Sigma_{n,d',n-1}
\cap H_0 \cap H_1 \cap \cdots \cap H_v$ in place of
$\Sigma_{n,e,n-1} \cap H_0,$
  where  $H_1$, \dots , $H_v$  are the branches of a hypersurface
in ${\Bbb P}^{N}$  corresponding to the virtually non-existent nodes.
 In fact, the singular curves form a hypersurface  ${\overline
V}_{n,1,0}\subset
{\Bbb P}^{N}$,  and each node of a curve of degree  $n$  determines
a unique branch of  ${\overline V}_{n,1,0}$ through that curve.  By
Proposition
1.2, the multiplicities  $b_j$'s  remain unchanged when
we intersect the ambient families with  $H_1$, \dots , $H_v.$

The same argument establishes (i) and (ii) as well.
\enddemo

 \remark {{\rm 3.10} Remark}  For  $2 \leq m \leq n$,  ${\overline
V}_{n,d,m}
 \cap\{a_{0n-m}= 0\}$
  contains components with general members of the form $C + l = C_1
+ \cdots + C_q + l,
\   l \nsubseteq C,$  and each  $C_r$  is irreducible, where
either  $q = 1$  or  $\deg(C_{r}) = 1$  for  $2 \leq r \leq q.$  Indeed,
if  $q \geq 2$  and
$m(C_1) \ne 0,$  then  $m(C_2) = 0$  by Proposition  3.9(iii).  We
suppose $\deg(C_2)
 \geq 2.$   The following surgeries establish the existence of the
required components.

 We will decrease the degree of  $C_2$  and increases the degree of
$C_1.$  First,
we degenerate  $C_2$  into a nodal curve  $C_2^{\prime} + L,$
 where  $L$  is a general line and  $C_2^{\prime}$  a sufficiently
general curve \cite{H}.  By the Principle of Degenerations,
$C_2^{\prime}+ L$
  has acquired several additional nodes.  We then consider $C_1 +
C_2^{\prime} + L$
and smooth several nodes of  $C_1 + C _2^{\prime}$  and  $C_1 + L.$

 Applying similar surgeries, one can easily obtain a list of
{\it all\/}  components of  $({\overline V}_{n,d,m} \cap\{a_{0n-m} =
0\})_{\text{red}}$
whose general members contain $l.$  We omit details.
\endremark

\head 4. Admissible schemes:  notation, definitions, and lemmas
\endhead
\subhead 4.1. The hyperplanes  $\{H_\sigma\}$ \endsubhead
 We consider the following sequence of  $n+1$  hyperplanes in
${\Bbb P}^{N}$:  $a_{0n} = 0$,
$a_{0n-1}= 0$, \dots , $a_{00} = 0.$  If  $G \in {\Bbb P}^{N }$ belongs
to the intersection
 of these hyperplanes, then
$$
f_{G}(X,Y,Z)  =  X(\sum a_{jk}X^{j-1}Y^{k}Z^{n-j-k})  =  X\{X(\cdots) +
\sum a_{1k}Y^{k}
Z^{n-1-k}\}.
$$
 Next,
we consider the following  $n$  hyperplanes in
${\Bbb P}^{N}$:  $a_{1n-1} = 0$,  $a_{1n-2} = 0$, \dots , $a_{10} = 0,$
etc.
We obtain  $\kappa$  hyperplanes in  ${\Bbb P}^{N}$, where
$\kappa = n + 1 + n + \cdots + 4.$
If    $G \in {\Bbb P}^{N }$ belongs to the intersection of these
$\kappa$  hyperplanes,
then  $G = (n-2)l + D.$  We shall employ a sequence
$\{H_{\sigma}\}$ of  $\kappa + 3$
  hyperplanes in ${\Bbb P}^{N}$,  where
$$
H_1= \{a_{0n} = 0\},\,  H_2= \{a_{0n-1} = 0\},\ldots , H_{\kappa+3}=
\{a_{n-20} = 0\}.
$$

\subhead 4.2. Standard exact sequences \endsubhead
Let  $V\subset {\Bbb P}^{N}$  be an arbitrary projective scheme and
$H \subset
{\Bbb P}^{N}$  a hyperplane.  In the sequel, we denote by  $h$
 a form defining  $H.$  For a positive integer  $k,$  let  $A^{k}(h)
\subset {\Cal O}_{V}$
denote the annihilator sheaf of the ideal
sheaf  $(h^{k}) \subset {\Cal O}_{V}.$  We have a standard exact
sequence of sheaves
on  $V$:
$$
0@>>>{\Cal O}_{V^k}(r-k)@>\otimes h^k>>{\Cal O}_V(r)@>>>{\Cal
O}_{V_k}(r)@>>>0,\tag{$3_k$}
$$
 where  ${\Cal O}_{V_k} ={\Cal O}_{V}/(h^{k})$  and  ${\Cal O}_{V^k} =
{\Cal O}_{V}/A^{k}(h).$
  This sequence yields a cohomology
sequence
$$
H^i({\Cal O}_{V^k}(r-k))@>>>H^i({\Cal O}_V(r))@>>>H^i({\Cal
O}_{V_k}(r)).\tag{$4_k$}
$$
 Therefore $H^i({\Cal O}_V(r)) = 0$ \,if \,$H^i({\Cal O}_{V^k}(r-
k))=H^i({\Cal O}_{V_k}(r))=0.$

Given a scheme  $V,$  we denote by  $\min \sup (V)$ the set of
{\it minimal associated points\/}  of  $V.$

\proclaim{4.3. Lemma} With the above notation, let  $  s = s(V,H)$ be
the smallest
integer such that $A^{s}(h) = A^{s+1}(h) = ...\  .$  Then  $V\backslash
V_{s} = V^{s }\backslash
V_{s}$ and
$$
\min \sup (V^s) = \min \sup (V) \backslash \sup(V_s)
\subset \min \sup (V^k),  \quad k \in {\Bbb Z}_{+}.
$$
\endproclaim

\demo{Proof} The problem is local and we may restrict everything
to  $\Cal O_{v,V};$  we denote the restrictions by  [$\cdot$].
 Clearly  $\min \sup (V^s)\backslash \sup(V_{s}) = \min \sup
(V)\backslash \sup(V_{s}).$
 It remains to show that  $\min \sup (V^s) \subset \sup(V)
\backslash \sup(V_{s}).$
  Let  $[A^{s}(h)] = \cap_{i }\, J_{i}$  be a minimal primary
decomposition.  Suppose
 $\text{rad} (J_{1}) \in [\min \sup (V^s)]$  and  $[h] \in
\text{rad}(J_{1}).$  Take  $a \in
\cap_{i\geq 2 }\,J_{i}\backslash J_{1}.$  Then  $a[h]^{e } \in
[A^{s}(h)]$  for
 $e>> 0.$  Hence   $a^{ }\in [A^{s}(h)],$  a contradiction.
\enddemo

\definition{4.4. Definition} A closed irreducible family of curves
of degree  $n$  is
{\it maximal\/} if its general members are of the form  $D + (n-w)l,$
where  $D$
is either a general member of a component $W \subset U_{m}(w,g)$
or  $D = 0.$
\enddefinition

\definition{4.5. Definition} We consider ${\overline V}_{n,d,m}$  with
$m \leq n.$  A closed
subscheme  $W \subset {\Bbb P}^{N}$  of dimension $\ge 2$  is $
{\overline V}_{n,d,m}$-{\it admissible\/} (or simply {\it
admissible})  if either
$W ={\overline V}_{n,d,m },$  or
$W = V^k$  or  $V_k$ in ($3_{k}$), where $V$  is admissible,
$V_{\text{red}} \nsubseteq H =
H_{\sigma}$ for the smallest possible  $\sigma = \sigma (V),$  and
$k$  is the smallest
integer such that   $\min \sup (V^k) = \min \sup (V) \backslash
\sup(V_k).$
\enddefinition

\remark {{\rm  4.6.} Remark} Let  $V$  be an admissible scheme and
$V_{\text{red}}
 \nsubseteq H = H_{\sigma}$  for the smallest possible  $\sigma .$
For every positive
integer  $s$  in ($3_s$), each component of  $(V_s)_{\text{red}}$ is a
maximal family by
Proposition 1.2(a) and Theorem 2.2.
 So, by Lemma 4.3, the integer  $k \leq s(V,H)$  and each component
of  $(V^{k})_{\text{red}}$  is a maximal family.
\endremark

Admissible schemes can be generically non-reduced as the following
example illustrates;  see also Section 3.

\example{4.7. Example} The family  ${\overline V}_{4,2,2}$  lies in
${\Bbb P}^{12}\ (\subset
{\Bbb P} ^{14}).$  The scheme  ${\overline V}_{4,2,2}\cap H_3$
contains an irreducible
 subscheme  $W$  whose general member is a curve of the form  $D
= C +l,$  where  $C$
  is a general cubic.  Exactly three transversal branches of
${\overline V}_{4,2,2}$  are
 passing through  $D,$ one for each pair of the nodes of  $D.$  By a
simple
dimension count  $\dim{\bold T}_{D}({\overline V}_{4,2,2}) = 12.$
Hence ${\overline V}_{4,2,2}\cap H_3$  is non-reduced at  D,  and
$\dim{\bold T}_{D}(W)
 = 11.$
\endexample

\definition{4.8. Definition} Set  $H_0= {\overline V}_{n,0,0}= {\Bbb
P}^{N}.$  Let  $V$
  be a  ${\overline V}_{n,d,m}$-admissible subscheme,
and  $K$  a component of  $V.$  We define a subset  $I(K) = I(K,V)
\subset {\Bbb Z}$  as follows:
$$
 \qquad i \in I(K) \quad  \iff \quad  K_{\text{red}}\subset H_i,  \
V\cap
H_0 \cap \cdots \cap H_{i-1} \nsubseteq H_i.
$$
\enddefinition

\proclaim{4.9. Lemma}  $\bold G({\overline V}_{n,d,m}) = {\overline
V}_{n,d,m}$  and
$\bold G(H_{\sigma}) = H_{\sigma}$ for  $\sigma = 1$, \dots ,
$\kappa+3;$ moreover,  $
\bold G(W) = W$ for every admissible $W.$
\endproclaim

\demo{Proof}  Since  $\bold G \subset PGL(2) \subset PGL(N)$,  we
get
$\bold G({\overline V}_{n,d,0}) = {\overline V}_{n,d,0}.$  By looking
at equations of
curves (see Section 1), we conclude that $\bold G(H_{\sigma})  =
H_{\sigma},$ hence
$\bold G({\overline V}_{n,d,p}) = {\overline V}_{n,d,p}$  if  $p \leq
n.$  Consider ($3_{k}$)
with an admissible  $V$  such that $\bold G(V)  = V$  and  $H =
H_{\sigma}.$  Since
$\bold G(H) = H,$ we get
$\bold G(A^{k}(h)) = A^{k}(h)$  for every  $k.$  Thus  $\bold G(W) =
W$  for every
admissible  $W.$
\enddemo

\subhead {4.10. Notation} \endsubhead  Let  $V$  be an admissible
scheme, and  $K$  an
irreducible subvariety of  $V$  with a general member
$$
C + \mu l = C_{1} + \cdots + C_{q} + \mu l,  \quad l \nsubseteq C,
$$
 where each  $C_{r}$  is irreducible  ($1 \leq r \leq q$).  We say that
K  is a subvariety of
{\it level\/}  $\mu ,$  {\it contact\/}  $m = m(C),$  and {\it genus\/}
$g = g(C).$
  Consider ($3_{k}$) as in Definition 4.5.  A component  $K^{\prime}$
of
 $V_{k}$  is said to be
{\it new\/}  if $ K_{\text{red}}^{\prime}$ is not a component of
$V_{\text{red}}.$

\definition{4.11. Definition}  a) With the above notation, a component
$K$  of  $V$
is said to be  $\delta_{\bold p}$-{\it nice\/}  if  $\delta_{\bold
p}(C_{r}) = 0 \  (1
\leq r \leq q).$

 b) A  $\delta_{\bold p}$-nice component  $R$  of a  ${\overline
V}_{n,d,m}$-admissible
scheme  $W$  is said to be  {\it nice\/}  if either
$R = W ={\overline V}_{n,d,m},$  or $W = V_{k}$  or  $V^{k }$ in
($3_{k}$)
of Definition 4.5 and  $R$  is a component of  $W$  coming from a
nice component  $M$  of  $V,$  that is,  $R_{\text{red}}$ is a
component
of  $(M \cap H_{\sigma (V)})_{\text{red}}$  or  $R_{\text{red}}=
M_{\text{red}}.$

 c) If  $K$  is a nice  ($\delta_{\bold p}$-nice) component of  $V,$
then  $K_{\text{red}}$  is said
to be a {\it nice\/}   ($\delta_{\bold p}$-{\it nice}) component of
$V_{\text{red}}.$
\enddefinition

\proclaim{4.12 Lemma} With the above notation, we have:

{\rm a)} If an admissible scheme  $V$ contains a component of
dimension  $\nu,$ level
 $\mu ,$ contact $m,$ and genus $g$ whose general members contain
$u$ general
lines through $\bold p$ and  $v$ general lines, then it contains an
irreducible subvariety
with the same data of a nice component of $V.$

{\rm b)} One and only one of the following conditions holds:

{\rm i)} all nice components of  $V_k$ are new and $K_{\text{red}}
\nsubseteq H_{\sigma(V)}$
 for every nice component  $K$ of $V;$

{\rm ii)} a nice component  $K^{\prime}$  of  $V_k$ is new if and
only if
$K_{\text{red}}^{\prime} \nsubseteq H_{\sigma(V_k)};$ in addition,
$V$ and $V_k$
 have nice components of different levels, and all nice
components of $V_{k }$ of smaller level are new;  or

{\rm iii)} $V_k$ has no new nice components and all its nice
components have the same level.

{\rm c)} If  $V$ has two nice components of level $\mu $   and
contacts $m$ and  $m_{2},$
 respectively, then it has a nice component of level  $\mu$ and
contact  $e$
for every  $e \ (m_{1} \leq e \le m_{2}).$  Moreover $m_{1} - g_{1} =
m_{2}
 - g_{2},$   where $g_{i}$   is the genus of the component of contact
$m_{i} \ (i = 1,2).$
\endproclaim

\demo{ Proof} The lemma is trivial if  $V={\overline V}_{n,d,m}.$
We proceed
by induction on the size of an admissible scheme.  Assume the
lemma for all admissible  $W \supseteq V.$  Then  $V^k$ satisfies
the lemma.

Now we turn to  $V_k.$  To begin with, we make the following
two remarks.

First, consider  ${\overline V}_{n,d,m}$  and  ${\overline
V}_{n,d,m+1}$
with  $2 \leq m \leq n - 1.$  Let  $U^{\prime}$
denote a union of the components of  $({\overline V}_{n,d,m}
\cap\{a_{0n-m}=
0\})_{\text{red}}$  that lie in  $\{a_{0n-m-1}= 0\}.$   By Proposition
1.2(a),
the general members of  $U^{\prime}$  contain  $l.$  Furthermore,
every component
of  $({\overline V}_{n,d,m+1} \cap\{a_{0n-m-1}=
0\})_{\text{red}},$  whose general members
contain $l,$  lies in  $U^{\prime}$  by Proposition 3.9.

Now, consider a component  $W^{\prime}$  of  ${\overline
V}_{n,d,m}\cap\{a_{0n-m}
 = 0\},\   m \leq n,$  whose general members contain $l.$  We
assume  $W^{\prime}$
  has contact  $m^{\prime} < \min \{d-n+m, m-2\};$  so
$W_{\text{red}}^{\prime}\subset
\{a_{1n-m^{\prime}} = 0\}.$  Let  $U^{\prime\prime}$  denote a
union of the components
of $({\overline V}_{n,d,m} \cap\{a_{0n-m}= 0\})_{\text{red}}$ that
lie in  $\{a_
{1n-m^{\prime}-1} = 0\}.$  Then  $(W^{\prime} \cap\{a_{1n-
m^{\prime}-1}= 0\})
_{\text{red}} \subset U^{\prime\prime}.$   In view
of Proposition 3.9, this follows at once by induction on  $d,$  the case
$d =
 n - m + m^{\prime} + 1$  being trivial.

{\it Completion of the proof.}\   First, we assume that
$K_{\text{red}}\nsubseteq
H_{\sigma(V)}$ for every nice component  $K$  of  $V.$
 Then all nice components of  $V_{k}$  are new, and  $V_k$ satisfies
(a) by induction hypothesis and Proposition 3.9.  Note that if
an irreducible curve   $E$  degenerates into  $C +l$  as in Theorem
2.2, the curve  $C$
  cannot contain lines through  $\bold p$  by Theorem 2.2 and the
Principle of
Degenerations, although it may contain general lines provided  $E$
has sufficiently
small genus (compare Remark 3.10).  Furthermore  $m(C) \leq m(E) -
2$  by Remark 2.4.
The assertion (c) also follows from Proposition
3.9 and the formula  $m^{\prime} - g^{\prime} = m - g - 2$  of
Theorem 2.2.

 Now, we assume that  $V$  contains nice components of different
levels.  Then, by (i) or (ii) for  $V$  (in place of  $V_{k}),$
 all nice components of  $V$  of smaller level are new and
$K_{\text{red}}\nsubseteq
 H_{\sigma(V)}$ for every nice component  $K$  of  $V$
of smaller level.   By the remarks,  $V_k$ satisfies (ii)
or (iii);  $V_k$ satisfies (iii) if and only if all nice components
of  $V_{k}$  have the same level.  Furthermore,  $V_{k}$  satisfies
(a) and (c) as in the preceding case.

 Finally, we assume that all nice components of  $V$  have the
same level and  $V$  has nice components of different contacts.
 By the remarks,  $V_k$  has no new nice components so it satisfies
(iii).  This proves the lemma.
\enddemo

Now, we come to the following key lemma which will enable us to
prove Theorem 5.2
(Vanishing Theorem).

\proclaim{4.13. Lemma} Let $W$   be an arbitrary  ${\overline
V}_{n,d,m}$-admissible
scheme, and $g \in W$  a general point. Let $W(g) \subset W$ be a
nice component whose
 general point is $g.$ Then
$$
      \dim{\bold T}_{g}(W(g))  \geq  \dim W(g) + \#I(W(g)). \tag{5}
$$
Furthermore, a basis of ${\bold T}_{g}(W(g)_{\text{red}})$   together
with  $\#I(W(g))$
elements of ${\bold T}_{g}(W(g))$ corresponding to the hyperplanes
$H_{i}\  (i \in I(W(g)))$
 form a linearly independent subset of ${\bold T}_{g}(W(g)).$
\endproclaim

\demo{Proof}  The lemma is trivial if  $W ={\overline V}_{n,d,m}.$
We proceed by induction
 on the size of  $W.$  Consider (3$_{s}$) as in Definition 4.5 (with $s =
k).$  Assume
the lemma for all admissible  $W\supseteq V.$  Then  $V^{s}$
satisfies the lemma.

Now we turn to  $V_s.$  Let  $V_{s}(g_s)$  be a nice component
of  $V_{s}$ with a general point  $g_{s}.$  Then  $V_{s}(g_s)$
 is coming from a nice component of  $V$  denoted by  $V(g_s).$

First, we consider the case when  $V_{s}(g_s)_{\text{red}}$  is
a component of  $V_{\text{red}},$  that is,  $V_{s}(g_{s})_{\text{red}}
\subset H_{\sigma}.$
  If  $s = 1$  then  $\#I(V_{s}(g_{s})) \leq \#I(V(g_{s})) - 1,$  and the
lemma follows by
induction.  But if  $s \geq 2,$  then  $\dim{\bold T}_{g_s}(V_{s}(g_s))
=
 \dim{\bold T}_{g_s}(V(g_s))$  and the lemma follows.

Next, we assume that  $g_{s}$  is  {\it not\/}  a general point of  $V$.
Let  $D_{s}$
denote the curve corresponding to  $g_{s}.$  Let  $D$  be a general
point
of  $V(g_s)$;  $D$  degenerates into  $D_{s}.$  We have
$$
 D = C + \mu l,  \qquad  D_{s} = C_{s} + \mu_sl,   \qquad (l \nsubseteq
C,\ C_{s}).
$$
By (5) for $V(g_{s}),$  we get
$$
 \dim{\bold T}_{D_s}(V(g_{s}))  \geq  \dim V(g_{s}) + \#I(V(g_s)).
$$
If  $\mu  =\mu _{s},$  then we consider two cases,  $s = 1$  and  $s
\geq 2,$  and deduce the lemma by induction.  For instance, if
 $s = 1$  then  $\#I(V_s(g_s)) \leq \#I(V(g_s)),$  since
 $V(g_s)_{\text{red}} \subset H_{i}$  for  $i \leq \sigma - 1$ (see
Definition 4.5).

Throughout the rest of the proof, we assume that  $\mu \ne \mu
_{s},$
 that is,  $\mu _{s} = \mu + 1.$  Let
$$
 C = E + F,    \qquad C_{s} = E_{s} + F,
$$
where
 $E$  is an irreducible curve that degenerates into  $E_{s} + l.$  We
claim:
$$
 \dim{\bold T}_{D_s}(V(g_s))  \geq  \dim{\bold T}_{D}(V(g_s)) + e +
\deg(E) - m(E),
\tag{6}
$$
where $e$  is the number of nodes of  $E$  approaching $\bold p$
 as  $E$  tends to  $E_s + l.$  In the proof of the claim, we can assume,
without
loss of generality, that  $V(g_s) = V(g_s)_{\text{red}}.$  Indeed,
$\dim{\bold T}_{v}(V(g_s))$  is an upper semicontinuous function
on the set of points of  $V(g_s)$  and nilpotents can only
"improve" the inequality.  Assuming  $C$  is reducible, we consider
the product
$$
\Pi= {\Bbb P}^{N_E} \times {\Bbb P}^{N_F}, \quad  N_{E} =
\deg(E)(\deg(E) + 3)/2,\
N_{F}= \deg(F)(\deg(F) + 3)/2,
$$
 and
a natural morphism  $\eta :\  \Pi \rightarrow  {\Bbb P}^{M},$  where
$M = \deg(C)
(\deg(C)+3)/2.$  (If  $C = E$ we set  $\Pi = {\Bbb P}^{N_E}.)$  Then  $E
+ F = \eta (E \times F),$
  and we can apply Proposition 3.9 to  ${\overline V}_{n,d,m} \subset
{\Bbb P}^{N_E},$
  where $ n = \deg(E)$,  $m = m(E),$  and $d$  is the number of nodes
of  $E.$
  This proves the claim.

 Now we will verify the following inequality:
$$
   \#I(V(g_s)) + e + \deg(E) - m(E) + 1 \ge \cases \#I(V_s(g_s))+1, &
s=1\\
\#I(V_s(g_s)), & s\ge 2.\endcases \tag{7}
$$
 Let  $i$  be the minimal integer in  $I(V_s(g_{s})).$  Then either
 $i \in I(V(g_s))$  or we get the vanishing of the corresponding
coefficient in the equation of  $E$  (see the proof of Proposition
3.9).  We then take the next integer in  $I(V_s(g_s))$
 and repeat the argument, etc.

 Finally, combining (6) with (5) for  $V(g_s)$  and (7), we
get (5) for  $V_s(g_s).$  The above discussion also proves
the last assertion of the lemma in case  $\mu_{s} =\mu + 1.$
\enddemo

\remark {{\rm 4.14.} Remark} Let $W$ be a component of
$U_m(n,g).$
One can now describe {\it all\/} components of  $W \cap \{a_{0n-m}
= 0  \}.$
Further, utilizing standard exact sequences it is not difficult to
calculate the
Hilbert polynomial of ${\overline V}_{n,d,m} \subset {\Bbb P}^{N-
m}$ as well as some other
interesting schemes of curves.
\endremark

\head{5. A vanishing theorem for admissible schemes}\endhead

\proclaim{5.1. Lemma} We consider $(3_{1})$  such that $V$  is
admissible, $V^{1}
\simeq V,\,  H^0(V_{1},{\Cal O}(-r))\mathbreak = 0$ for $r >> 0$,  and
$H^1(V_1,{\Cal O}(r)) = 0$
 for all $r \in {\Bbb Z}.$  Then $H^0 (V ,{\Cal O}(-r )) = 0$ for  $ r >>
0,$ and $
H^1 (V ,{\Cal O}(r )) = 0$  for all  $r \in {\Bbb Z}.$
\endproclaim
\demo{Proof}  Recall first the following classical lemma of
Enriques - Severi - Zariski  \cite{G, Exp XII, Corollary 1.4}: {\it for an
arbitrary
projective scheme\/} $V \in {\Bbb P}^{N},\  \text{depth}({\Cal
O}_{v,V}) \geq 2$
{\it for every point\/} $v \in V$ {\it if and only if\/} $H^{i}(V,{\Cal
O}(-r)) = 0$
{\it for\/} $r >> 0$
{\it and\/}  $i = 0,1.$

 It follows that  $\text{depth}({\Cal O}_{C,V_1 }) \geq 2$  for every
point  $C \in V_{1},$
  hence               $\text{depth}({\Cal O}_{C,V}) \geq 2$  for every
point  $C \in V_{1}
\subset V.$

 Let  $C \in V$  be an arbitrary point.  We will show that  $[1:0:0]
\notin \gamma (C)$  for a suitable  $\gamma \in {\bold G }\subset
PGL(N).$
Assuming  $[1:0:0] \in C,$  take a point  $[1:a:b] \notin C$  with  $a, b
\in {\Bbb Z}_{+}.$
It is easy to find  $\gamma \in {\bold G } \subset PGL(2)$  such that
$\gamma ([1:0:0])
 = [1:a:b].$  Then  $[1:0:0] \notin \gamma^{-1}(C).$

 Now we assume that  $[1:0:0] \notin C.$  Consider the points
$(\varphi
_{tt}\cdot \phi_{tt})(C)$  for  $t \in {\Bbb C}^*.$  As  $t$  goes to  $0$,
$
(\varphi
_{tt}\cdot \phi_{tt})(C)$  tends to a point  $nl$  of the subscheme
$V_1 \subset V.$
  Since  $\text{depth}({\Cal O}_{C,V})$  is an upper semicontinuous
function
of  $C, \ \text{depth}({\Cal O}_{C,V}) \geq 2$  for every point  $C \in
V.$
 Thus, Lemma 5.1 follows from the Enriques - Severi - Zariski
lemma and ($4_k$).
\enddemo

\proclaim{5.2. Theorem (Vanishing Theorem)} Let $V \subset {\Bbb
P}^{N}$   be a
${\overline V}_{n,d,m}$-admissible subscheme.  Then  $H^1 (V ,{\Cal
O}(r )) = 0$ for all $
r \in {\Bbb Z},$ and $H^0(V ,{\Cal O}(-r)) = 0$ for  $r >> 0.$
\endproclaim

\demo{Proof}{\it Case\/}:  $\dim V = 2.$  We consider six closed
irreducible families
in  ${\Bbb P} ^{5}\ (\subset {\Bbb P}^{N}),$  whose general members
are of the
form  $D + (n-2)l$ with  $l \nsubseteq D,$  by specifying  $D$:

${\Cal F}_{1}$: $D$  is a general quadric; \ \   ${\Cal F}_{2}$: $D$  is a
general quadric with
$(D\cdot l)_{\bold p} = 1;$

${\Cal F}_{3}$: $D$  is a general quadric with  $(D\cdot l)_{\bold p} =
2;$

${\Cal F}_{4}$: $D$  is a union of two general lines;

${\Cal F}_{5}$: $D$  is a union of two general lines with $(D\cdot
l)_{\bold p} = 1;$

${\Cal F}_{6}$: $D$  is a union of two general lines with $(D\cdot
l)_{\bold p} = 2.$
\smallskip
We also consider the closed irreducible family  ${\Cal F}_7 \subset
{\Bbb P}^{5}$  whose
general member is of the form  $D + (n-1)l,$  where  $D$  is a general
line.  There are only
four maximal families of dimension at most  $2,$  namely:  ${\Cal
F}_{6}$,  ${\Cal F}_{7}$,
$\{(n-1)l + l' \  |\  {\bold p} \in  l^{\prime} \}$,  and  $\{nl \}.$  The
latter two families
are contained in any admissible scheme.  We get
$$
 V_{\text{red}} \subset {\Cal F}_{6 } \cup {\Cal F}_{7} \subset
\bigcap_{1\leq \sigma\leq
\kappa +2}  H_{\sigma}= {\Bbb P}^{3} \subset {\Bbb P}^{5} \subset
{\Bbb P}^{N},
\qquad {\Cal F}_{7 }= {\Bbb P}^{3 } \cap H_{\kappa+3}.
$$
Moreover $H^1 (V_{\text{red}},{\Cal O}(r ))  = 0$  for all $r.$  Thus we
get the vanishing
of the cohomology groups provided  $V \subset {\Bbb P}^{3}.$  If
$V_{\text{red}}
={\Cal F}_{6 }\cup {\Cal F}_{7},$  we take  $H = H_{\kappa+3}$  and
consider ($3_{k}$)
for an appropriate
integer  $k.$  Then  $(V_k)_{\text{red}}={\Cal F}_{7}$  and
$(V^k)_{\text{red}}=
{\Cal F}_{6}.$

 Now, we assume that  $V$  is irreducible, that is,
$V_{\text{red}}={\Cal F}_{6 }$ or
${\Cal F}_{7}.$  First we take the hyperplane  $H = H_{1}.$  Let  $h =
0$  be an equation
of  $H.$  Consider the following exact sequence similar to $(3_{1}):$
 $$
0@>>>K@>\otimes h>>{\Cal O}_V(r)@>>>{\Cal O}_{V_1}(r)@>>>0,
$$
where $K = 0$  if  $h{\bigm |}_{V}= 0,$  and  $K = {\Cal O}_{V^1}(r-1)$
otherwise.
In the latter case,  $\dim  V^{1}   = 2$  by Lemma 4.13.  Indeed, the
restriction
of  $h$  to the local ring of  $V$  at its general point is a nontrivial
element.  Hence the stalk of  $A^{1}(h)$  at this point is a nontrivial
ideal
in the corresponding local ring.  This ideal is nilpotent, so $\dim
V^{1}  = 2.$  Next, we
consider  $V_{1}$ (or  $V^{1}$)  in place of  $V$  and take the
hyperplane
 $H = H_2$ (or  $H = H_1$),  etc.  After a finite number of steps, we
obtain subschemes in
${\Bbb P} ^{3}.$

{\it  General case.} Assume the theorem for all smaller admissible
schemes.  If  $\dim V \geq 3,$  we take the hyperplane
$H_{\sigma}$
 with the smallest possible  ${\sigma}$  such that  $V_{\text{red}}
\nsubseteq H_{\sigma}.$
  For an appropriate  $k,$  we consider ($3_{k}$) with
 $H = H_{\sigma}$  and get two admissible schemes,  $V_{k}$  and
$V^{k}.$  Since  $V_{k} \ne  V,$  we can apply Lemma 5.1 if
$V^{k}\simeq V.$
This proves the theorem.
\enddemo

\enddocument